\documentstyle[12pt]{article}
 \def\Bbb{\bf}

 \catcode`@=11
 \def\@date{January 15, 1995}
 \catcode`@=12
\newcommand{\sgn}{\mathop{\rm sgn}}
\newcommand{\Row}{\mathop{\rm Row}}
\newcommand{\Col}{\mathop{\rm Col}}

\newcommand{\IM}{\mathop{\rm Im}}
\newcommand{\Sym}{\mbox{\rm Sym}\mbox{}}

\newcommand{\Gr}{\mbox{\rm Gr}}
\newcommand{\id}{\mathop{\rm id}}

\newcommand{\ind}{\mbox{\rm ind}}

\newcommand{\res}{\mbox{\rm res}}

\newcommand{\wt}{ {\mathop{\rm wt}}\mbox{} }
\newcommand{\tr}{\mathop{\rm tr}}

\newcommand{\Hom}{\mbox{\rm Hom}}
\newcommand{\Char}{ {\mbox{\rm char}} }
\newcommand{\Span}{ {\mbox{\rm Span}} }
\newcommand{\diag}{\mbox{\rm diag}}
\newcommand{\Div}{\mathop{\rm div}}
\newtheorem{thm}{Theorem}
\newtheorem{lem}[thm]{Lemma}
\newtheorem{cor}[thm]{Corollary}
\newtheorem{prop}[thm]{Proposition}

\newcommand{\eqdef}{\stackrel{\rm def}{=}}
\newcommand{\CC}{{\Bbb C}}

\newcommand{\ZZ}{{\Bbb Z}}
\newcommand{\NN}{{\Bbb N}}
\newcommand{\PP}{{\Bbb P}}
\newcommand{\al}{{\alpha}}

\newcommand{\alD}{{\alpha_D}}
\newcommand{\beD}{{\beta_D}}

\newcommand{\lam}{ \lambda }
\def\SS{ { \cal S\hspace{-.2em}S } }
\newcommand{\FM}{ { \cal F\hspace{-.2em}M } }
\newcommand{\FF}{ { \cal F } }
\newcommand{\FD}{ \FF_D }
\newcommand{\FSS}{{ \FF_{3,n}^{\SS} }}
\newcommand{\FFM}{{ \FF_{3,n}^{\FM} }}
\newcommand{\FFMt}{{ \FF_{3,3}^{\FM} }}
\newcommand{\Ft}{{ \FF_{3,n} }}
\newcommand{\Ftt}{{ \FF_{3,3} }}
\newcommand{\II}{ { \cal I } }

\newcommand{\LL}{ {\cal L} }
\newcommand{\LD}{ \LL_D }
\newcommand{\OO}{ {\cal O} }

\newcommand{\Pijk}{ \PP^1_{ijk} }
\newcommand{\Bl}{\mbox{\rm Bl}}
\newcommand{\ReRo}{\mbox{RR}}
\newcommand{\vv}{\mbox{\bf v}}
\newcommand{\down}[1]{\downarrow\rlap{\raise.3ex\hbox{$\!\scriptstyle{#1}$}}
}

\title{The Space of Triangles, Vanishing Theorems, and Combinatorics}
\author{Wilberd van der Kallen and Peter Magyar}

\begin{document}

\maketitle

\begin{center} {\bf Abstract} \end{center}

{\small
\noindent
We consider compactifications of
$(\PP^n)^3 \setminus \bigcup \Delta_{ij}$,
the space of triples
of distinct points in projective space.
One such space
is a singular variety of configurations of points and lines;
another is the smooth compactification of Fulton and MacPherson;
and a third is the triangle space of Schubert and Semple.

We compute the sections of line bundles on these spaces,
and show that they are equal as $GL(n)$ representations to the
generalized Schur modules associated to ``bad''
generalized Young diagrams with
three rows (Borel-Weil theorem).
On the one hand, this yields Weyl-type character and dimension
formulas for the Schur modules;
on the other, a combinatorial picture of the space of sections.
Cohomology vanishing theorems play a key role in our analysis.
}
\\[1em]
{\large \bf Introduction}
\\[1em]
{}First, we present our main character
in three different guises.

{\bf Combinatorial:}  the generalized Young diagram
$$
D_3 = \begin{array}{ccccccc}
\Box &      &      &      & \Box & \Box  & \Box  \\
     & \Box &      & \Box &      & \Box  & \Box  \\
     &      & \Box & \Box & \Box &       & \Box
\end{array}
$$
whose columns correspond to the non-empty subsets of set of row-indices
$\{1,2,3\}$.
Also, we fix an integer $n \geq 3$.

{\bf Geometric:} the configuration space
$\Ft$, which is the variety of all 7-tuples
$(p_1, p_2, p_3, l_1, l_2, l_3, P)$,
where $p_i$ are points, $l_i$  lines, and $P$
a plane, all in $\PP^{n-1}$; subject to
certain inclusions:
the point $p_i$ must lie on each of the lines $l_j$ for
$j \neq i$, and all must be contained in the plane $P$.
This the
{\em singular space of triangles in} $\PP^{n-1}$,
a projective variety of dimension $3n-3$.
The seven entries correspond to the seven columns of
the diagram, and each entry is a space whose projective dimension
is one less than the number of boxes in the column.
The inclusions of spaces
correspond to horizontal inclusions of columns.
(If we replaced $D_3$ by a Young diagram,
the corresponding configuration
space would be a flag variety.)

{\bf Algebraic:}  The $GL(n)$-representation $S_{D_ ,n}$,
a module
defined by the Schur--Weyl construction applied to $D = D_3$.
That is, let $\CC^n$ be the defining
representation of $GL(n)$, and consider the tensor power
$(\CC^n)^{\otimes D}$ with one factor for each square of the diagram.
Then
$$
S_{D_3,n} = (\CC^n)^{\otimes D} \gamma_D,
$$
where $\gamma_D$ is a Young operation symmetrizing and anti-symmetrizing
tensors according to the rows and columns of $D = D_3$.
(If $D$ were
a Young diagram, $S_{D,n}$ would be a
classical Schur module, an irreducible representation of $GL(n)$.)
%\\[1em]

\

\

These objects have been extensively explored.
Combinatorists have examined
the representations $S_{D,n}$ as part of the theory of generalized
Schur modules and Young diagrams.
(See \cite{MP}, \cite{RS1}, \cite{RS2}, \cite{RS3}, \cite{W}.)
The geometric theory of
$\Ft$ and its desingularizations goes back to Schubert \cite{Schubert}
and Semple \cite{Semple}, and has been illuminated recently by Fulton,
MacPherson, and others (\cite{CF}, \cite{FM}, \cite{LeBarz},\cite{RobSpeis},
\cite{BD1}, \cite{BD2}).

In the current paper, we explore the relations between
the algebraic and geometric pictures.  We prove a Borel-Weil
theorem realizing $S_{D_3,n}$
(or the Schur module of any three-row diagram) as the
sections of a line bundle
over the triangle space, and we derive explicit character and dimension
formulas generalizing those of Weyl for the irreducibles.
Reversing the perspective, the Schur modules give a
combinatorial construction of the space of sections of
certain line bundles over the triangle spaces.

The key steps relating the geometry to the combinatorics
require the vanishing of certain
higher cohomology groups
and explicit computations on
spaces of sections.
{}For these, we make a detailed examination of the defining
equations, torus-fixed points, and birational maps
of the triangle spaces, in order to apply the standard
techniques for dealing with (almost) homogeneous varieties.
Indeed, we intend this
paper partly as a primer on Frobenius splitting,
natural desingularizations,
Lefschetz theorems,
and rational singularities,
as illustrated by our simple example.

Nevertheless, our character formulas for Schur modules
are stated in purely elementary
terms, and the interested reader can skip directly
to sections \ref{character} and \ref{dimension}.

In previous papers \cite{MaNW}, \cite{MaSchub}, \cite{MaFour}
we considered
the same problems for diagrams $D$ satisfying the ``northwest''
or ``strongly separated'' conditions,
for which the geometry of $\FD$ is particularly
simple.  $D_3$  is the smallest diagram which does not fall
into these classes, and thus needs a different
treatment.
In the case of northwest or strongly separated diagrams,
we may desingularize our varieties by
the usual Bott--Samelson and
Zelevinsky resolutions, but for the triangle space
we must use more
complicated special desingularizations which appear in the
literature. We hope that our methods will shed light on defining
non-singular spaces of tetrahedra and higher-dimensional simplices.
Even for the space of triangles, some of the combinatorial results
seem geometric in nature, but have no obvious geometric explanation.
See especially section \ref{virtual}.

\vspace{1em}

\noindent
{\small
{\bf Contents.}\ \ \
{\bf 1.} Definitions\ \ \
{\bf 2.} The Singular space \
{\bf 2.1} Borel-Weil theorem \
{\bf 2.2} Frobenius splitting \
{\bf 2.3} Proof of Thm \ref{vanishing} \
{\bf 2.4} Fixed points \
{\bf 3.} The Schubert--Semple space \ \ \
{\bf 4.} The Fulton--MacPherson space \
{\bf 4.1} Strata \
{\bf 4.2} Fixed lines \
{\bf 4.3} Rational singularities \ \ \
{\bf 5.} Fixed-point formulas \
{\bf 5.1} General theory \
{\bf 5.2} Character formula \
{\bf 5.3} Dimension formula \
{\bf 5.4} Virtual desingularization  \ \ \
{\bf 6} Appendix: Restriction and induction \ \ \
References
}

\vspace{2em}

\noindent
{\bf Acknowledgements.}
The second author wishes to thank
V. Reiner and M. Shimozono, who introduced him to
these questions; and W. Fulton, for helpful references.

\section{Definitions}
\label{Definitions}

A {\em diagram} is a finite subset of
$\NN \times \NN$.
Its elements $(i,j) \in D$ are called {\em squares}, and
the square $(i,j)$ is pictured in the
$i^{\mbox{th}}$ row and $j^{\mbox{th}}$ column.
We shall often think of $D$ as a sequence
$(C_1,C_2,\ldots,C_r)$ of columns $C_j \subset \NN$.
% The Young diagram corresponding to
% $\lam = (\lam_1 \geq \lam_2 \geq \cdots \geq \lam_n \geq 0)$
% is the
% set $\{(i,j) \mid 1\leq j\leq n,\  1\leq i\leq \lam_j \}$.

{}Fix once and for all an integer $n\geq 3$.
We denote the interval $[1,n] = \{1,2,\ldots,n\}$.
We shall always write $G = GL(n,\CC)$, $B = $ the subgroup
of upper triangular matrices, and $T =$ the subgroup of diagonal
matrices.  $\CC^n$ is the defining representation of $G$.

We will assume our diagrams have at most $n$ rows:
$D \subset [1,n] \times \NN$.
Let $\Sigma_D$ be the symmetric group permuting the squares of $D$,
and for any diagram $D$, let
$$
\Col(D) = \{\pi \in \Sigma_D \mid \pi(i,j) = (i',j) \ \exists
 i'\}
$$
be the group permuting the squares of $D$
within each column, and we define $\Row(D)$ similarly
for rows.
Define the idempotents
$\al_D$, $\beta_D$ in the group algebra
$\CC[\Sigma_D]$ by
$$
\alD = {1 \over |\Row D|} \sum_{\pi \in \Row D} \pi, \ \ \ \
\beD = {1 \over |\Col D|} \sum_{\pi \in \Col D} \sgn(\pi) \pi ,
$$
where $\sgn(\pi)$ is the sign of the permutation.

Now, $G$ acts diagonally on the left of
the $|D|$-fold tensor product $(\CC^n)^{\otimes D}$,
and $\Sigma_D$ acts on the right by permuting the tensor factors:
$$
g(x_{t_1}, x_{t_2}, \ldots ) \pi =
(g x_{\pi t_1}, g x_{\pi t_2}, \ldots )  .
$$
These actions commute.
Define the {\em Schur module}
$$
S_D \eqdef (\CC^n)^{\otimes D} \alD \beD \subset (\CC^n)^{\otimes
 D},
$$
a representation of $G$.
Interchanging two columns (or two rows) of the diagram
gives an isomorphic Schur module.
%  if we change the diagram by permuting the rows or the
%  columns (i.e., for some permutation $\pi: \NN \rightarrow \NN$,
%  changing $D = \{(i,j)\}$ to
%  $D' = \{(\pi(i), j) \mid (i,j) \in D\}$,
%  and similarly for columns).

Now,  let
$e_1, \ldots, e_n$ be the standard basis of $\CC^n$,
and for $C\subset [1,n]$, define the coordinate subspace
$E_C = \Span(e_i \mid i\in C) \in \Gr(|C|,\CC^n)$.
{}For a diagram $D = (C_1,\ldots,C_r)$,
let
$$
\Gr(D) = \Gr(|C_1|,\CC^n) \times \cdots \times \Gr(|C_r|,\CC^n)
 \ ,
$$
and define the {\em configuration variety} as the closure
in $\Gr(D)$ of the $GL(n)$-orbit
of a configuration of coordinate subspaces:
$$
\FD = \mbox{closure } \left[ {G\cdot (E_{C_1},\ldots,E_{C_r})}
 \right]
\subset  \Gr(D) \ .
$$
This is clearly an irreducible subvariety.
We may also define the  {\em inclusion variety}
$\II_D \subset \Gr(D)$
by:
$$
\II_D = \{ (V_1,\ldots,V_r) \in \Gr(D)
\mid C_i \subset C_j \Rightarrow V_i \subset V_j \} .
$$
We clearly have $\FD \subset \II_D$.

Consider the Pl\"ucker line bundle
$\OO(1) = \OO(1,\ldots, 1)$ on the
product of Grassmannians $\Gr(D)$.
We may define a line bundle $\LD$ on $\FD$ and $\II_D$ as the
restriction of $\OO(1)$ to the subvarieties.

In case $D$ is a Young diagram $\{(i,j)\mid 1\leq i \leq \lam_j\}$
for $\lam = (\lam_1 \geq \cdots \geq \lam_n \geq 0)$, then
$\FD = \II_D$ is a flag variety, and $\LD$ is the Borel-Weil line
bundle whose sections are (the dual of) the irreducible Schur
module $S_D = S_{\lam}$.

In this paper, we will consider diagrams
$D$ with at most three rows:
that is,
all squares $(i,j) \in D$ have $i = 1$, 2, or 3.
We may consider
the diagram $D_3$ of the introduction
as  universal:
it contains a column of each type, and
any three-rowed diagram can be specified
(up to the order of the columns)
by the multiplicity
$m_i \geq 0$, $i= 1,\ldots, 7$ of each column in $D_3$.

Note that the multiplicities $m_i$ will not affect the
varieties $\FD$ or $\II_D$ (provided all $m_i > 0$).
We will denote $\FD = \FF_{D_3} = \Ft$ and
$\II_D = \II_{D_3} = \II_{3,n}$,
both inside the product of Grassmannians:
$$
\Ft \subset \II_{3,n} \subset (\PP^{n-1})^3 \times
\Gr(1,\PP^{n-1})^3 \times \Gr(2,\PP^{n-1}).
$$
The variety $\Ft$ is the closure of the $GL(n)$-orbit of the
coordinate 2-simplex in $\PP^{n-1}$, and the inclusion variety is
$$
\II_{3,n} =
 \{\ (p_1, p_2, p_3, l_1, l_2, l_3, P) \ \mid \ p_i \subset l_j
\ \ \forall\  i\neq j, \ \ l_j \subset P \ \ \forall \ j \ \}.
$$
{}For a diagram $D$ defined by integers $m_i$,
$\LD$ is the restriction of $\OO(m_1,\ldots,m_7)$ on the product
 of
Grassmannians.  This is a (very) ample line bundle
exactly when $m_i > 0$ for all $i$.

\section{The Singular space}

\subsection{Borel-Weil theorem}

Our first aim is to prove
\begin{thm}
\label{vanishing}
% which contains
% each of the columns of $D_3$ at least once,
We have \\
(a) For any diagram $D$ with at most three rows, we have
$H^0(\Ft, \LD) = S_D^*$ and $H^i(\Ft, \LD) = 0$ for $i >0$. \\
(b) $\Ft = \II_{3,n}$ is a normal, irreducible variety,
and projectively normal with respect to $\LD$.
\end{thm}

The proof will occupy the following sections.  We start with the
 elementary
parts.
\\[.3em]
{\bf Claim:\ } $\II_{3,n}$ is irreducible of dimension
$3n - 3$, and hence equal to its subvariety $\Ft$.
\\[.3em]
It will suffice to prove this for $\II_{3,3}$, since there is an
obvious fiber bundle
$$
\begin{array}{ccc}
\II_{3,3} & \rightarrow & \II_{3,n} \\
& & \downarrow \\
& & \Gr(2,\PP^{n-1})
\end{array}
$$
given by mapping a triangle to the
projective plane in which
it lies.
(The fiber is irreducible if and only if the whole bundle is, and
$\dim \II_{3,n} = \dim \II_{3,3} + 3(n - 3)$.)

Now, any triangle in $\II_{3,3}$ can be deformed under the
$GL(3)$ action on $\PP^2$ so as to approach a maximally degenerate
configuration, in which all three vertices and edges are identical.
{}Furthermore, $GL(3)$ acts transitively on this stratum
of degenerate triangles, so it suffices to check the irreducibility
in a neighborhood of a single degenerate triangle such as
$(E_1,E_1,E_1,E_{12},E_{12},E_{12})$,
for which the entries form a flag of coordinate subspaces of $\PP^2$.
Nearby, we can give affine coordinates so that
a configuration in $(\PP^2)^3 \times \Gr(1,\PP^2)^3$ is represented
by
$$
\left[
\left( \begin{array}{c} 1 \\ a_1 \\ b_1 \end{array} \right),
\left( \begin{array}{c} 1 \\ a_2 \\ b_2 \end{array} \right),
\left( \begin{array}{c} 1 \\ a_3 \\ b_3 \end{array} \right),
\left( \begin{array}{c} c_1 \\ d_1 \\ 1 \end{array} \right),
\left( \begin{array}{c} c_2 \\ d_2 \\ 1 \end{array} \right),
\left( \begin{array}{c} c_3 \\ d_3 \\ 1 \end{array} \right)
\right]
$$
subject to the six incidence conditions:
$$
\begin{array}{c}
c_2 + a_1 d_2 + b_1  =  0 \\
c_3 + a_1 d_3 + b_1  =  0
\end{array} \ \ \
\hfill
\begin{array}{c}
c_1 + a_2 d_1 + b_2  =  0 \\
c_3 + a_2 d_3 + b_2  =  0
\end{array}
\hfill \ \ \
\begin{array}{c}
c_1 + a_3 d_1 + b_3  =  0 \\
c_2 + a_3 d_2 + b_3  =  0
\end{array}
$$
By eliminating, we can reduce this to
the seven variables $a_1$, $a_2$, $a_3$, $d_1$, $d_2$, $d_3$, $c_1$
subject to the single equation
$$
(a_1 - a_2)(d_2 - d_3) + (a_3 - a_2)(d_1 - d_2) = 0 .
$$
A linear change of variables turns this into the product of an
affine space and a quadric, an irreducible variety of dimension 6.
This is what we wanted to show.
Therefore $\II_{3,n}$ is irreducible and equal to $\Ft$.
\\[.5em]
{\bf Borel-Weil construction.}
\\[.5em]
Next, we recall from \cite{MaNW} the elementary construction connecting
the Schur module $S_D$ of a three-row diagram with the triangle space.
Let $V = \CC^n$ and $U = V^*$ its dual space.
By definition, $S_D$ is the image of the composite map
$$
S_D = \IM\left[ V^{\otimes D} \alD
\stackrel{\mbox{\small incl}}{\rightarrow}
 V^{\otimes D}
\stackrel{\beD}{\rightarrow}
 V^{\otimes D} \beD
\right].
$$
Taking dual spaces, we have
$$
S_D^* = \IM\left[ U^{\otimes D} \beD
\stackrel{\mbox{\small incl}}{\rightarrow}
 U^{\otimes D}
\stackrel{\alD}{\rightarrow}
 U^{\otimes D} \alD
\right].
$$

We translate this into geometric language as follows.
Consider the product space $(\PP^{n-1})^D$ as
all $|D|$-tuples of points
inscribed in the squares of $D$.
Define
$$
\phi : (\PP^{n-1})^D \rightarrow \Gr(D)
$$
to be the rational map taking a $c$-tuple of vectors
in a column $C$ to the space in $\Gr(c,\CC^n)$ which they span.
This is a rational map defined everywhere except a set of codimension
 $\geq
2$, so
it induces maps of locally free coherent sheaves as if it were regular.
% not true without "locally free". Condider example of open subset
% whose
% complement Z has codim 2. The ideal sheaf of Z is isomorphic with
% the
% structure sheaf, but the isomorphism does not extend.
Also define the {\em row multidiagonal}
$\Delta^D(\PP^{n-1})$, the locus in $(\PP^{n-1})^D$
 where all points in the
same row are equal: that is, for our $D = D_3$,
$$
\Delta^D(\PP^{n-1}) =
\left\{ \begin{array}{ccccccc}
p_1 &      &      &      & p_1  & p_1  & p_1  \\
     & p_2 &      & p_2  &      & p_2  & p_2  \\
     &     & p_3  & p_3  & p_3  &      & p_3
\end{array}
\right\}
\subset (\PP^{n-1})^D \ \ .
$$
The composite image of these maps is precisely the configuration
 variety:
$$
\FD = \mbox{closure } \IM \left[
\Delta^D(\PP^{n-1})
\stackrel{\mbox{\small incl}}{\rightarrow}
 (\PP^{n-1})^D
\stackrel{\phi}{\rightarrow} \Gr(D)
\right] .
$$
The above equation for the dual Schur module
now translates easily into
$$
\begin{array}{rcl}
S_D^* & = & \IM\left[
\begin{array}{rcl}
H^0(\Gr(D), \LD)
& \stackrel{\phi^*}{\rightarrow} &
H^0( (\PP^{n-1})^D, \OO(1) )
\\
& \stackrel{\mbox{\small rest}}{\rightarrow} &
H^0( \Delta^D(\PP^{n-1}), \mbox{rest}\,\OO(1) )
\end{array}
\right]
\\
& & \\
& = &
\IM \left[
H^0(\Gr(D), \LD) \stackrel{\mbox{\small rest}}{\rightarrow} H^0(\FD,
 \LD)
\right] \ ,
\end{array}
$$
where $\mbox{rest} = \mbox{incl}^*$.
This equation is true for an arbitrary diagram $D$.
If $\FD = \Ft$, the first equation in part (a) of the Theorem states
 that
 the above
restriction map is onto, so that the dual Schur module is just
equal to the sections of $\LD$ over the triangle space $\FD = \Ft$.

If $D$ does not contain each column of $D_3$, then $\FD \neq \Ft$,
 but there
is a surjective map $\Ft \rightarrow \FD$ given by forgetting the
 data
associated
to the columns which do not appear in $D$.
We have the commutative diagram
$$
\begin{array}{ccc}
H^0(\Gr(D_3), \LD) & \stackrel{\mbox{\small rest}}{\rightarrow} &
 H^0(\Ft,
 \LD) \\
\downarrow & & \downarrow \\
H^0(\Gr(D), \LD) & \stackrel{\mbox{\small rest}}{\rightarrow} & H^0(\FD,
 \LD)
\end{array}
$$
The vertical map between the sections over Grassmannians is clearly
 an
isomorphism.
Hence in this case, we must show the surjectivity
of the top restriction map, {\em and} the bijectivity of the second
vertical map.

Thus, the first equation of part (a) reduces in general to
\begin{lem}
The natural map
$
H^0(\Gr(D), \LD) {\rightarrow} H^0(\FD, \LD)
$
is surjective, and the natural map
$
H^0(\Ft, \LD) \rightarrow H^0(\FD, \LD)
$
is bijective.
\end{lem}

To prove these facts and the rest of the Theorem, we will need more
sophisticated techniques.

\subsection{Frobenius splitting}

The theory of Frobenius splittings invented by Mehta,
Ramanan, and Ramanathan
(\cite{MR}, \cite{RR}, \cite{R1},
% \cite{R2},
\cite{vdK})
is a characteristic-$p$ technique for proving surjectivity
and vanishing results about coherent sheaves,
even in characteristic zero.
It is highly practical for dealing with homogeneous varieties
because one can work over the integers in a characteristic-free
way, and never consider the special features of characteristic-$p$
geometry.  In fact, the method reduces to classical questions
about defining equations and canonical divisors of varieties.
Most of the theorem of the last section will follow immediately from
 knowing
that the pair $\Ft \subset \Gr(D_3)$ is ``compatibly Frobenius split''
in any characteristic.

Given two algebraic varieties $Y \subset
 X$ defined
over an algebraically closed
field $F$ of characteristic $p > 0$,
with $Y$ a closed subvariety of $X$,
we say that the pair $Y \subset X$ is
{\em compatibly Frobenius split} if: \\
(i) the $p^{th}$ power map $F: \OO_X \rightarrow F_* \OO_X$
has a splitting, i.e. an $\OO_X$-module morphism
$\phi: F_* \OO_X \rightarrow \OO_X$ such that $\phi F$ is the
identity; and \\
(ii) we have $\phi(F_* I) = I$, where
$I$ is the ideal sheaf of $Y$.

Because $\LL\otimes F_* \OO_X=F_* \LL^p$ for any line bundle $\LL$,
a Frobenius splitting of $Y$ allows one to embed the
cohomology of an ample bundle on $Y$ into the cohomology
of its powers: $H^i(Y,\LL) \subset H^i(Y,\LL^{p^d})$
for all $d \geq 0$.  Since the right-hand side becomes
zero for large $d$ by Serre vanishing,
the $H^i(Y,\LL)$ itself must be zero ($i>0$).
If one can show this vanishing for reductions modulo
$p$ for all (or infinitely many) $p$,
then semi-continuity implies $H^i(X(\CC), \LL) = 0$ as well.
In fact, Mehta and Ramanathan prove the following
\begin{prop}\label{splitsur}
Let $X$ be a projective variety, $Y$ a closed subvariety,
and $\LL$ an ample line bundle on $X$.
If $Y \subset X$ is compatibly split, then
$H^i(Y,\LL) = 0$ for all $i >0$,
and the restriction map
$H^0(X,\LL) \rightarrow H^0(Y,\LL)$ is surjective.

{}Furthermore, if $Y$ and $X$ are defined and projective over $\ZZ$
(and hence over any field), and they are compatibly split
over any field of positive characteristic,
then the above vanishing and surjectivity
statements also hold for all fields of characteristic zero.
\end{prop}

{}Frobenius splitting is also sufficient to establish
the normality of our varieties.
The main theorem of
Mehta and Srinivas \cite{MS} states that if $Y$ is a Frobenius-split
variety possessing a desingularization with connected fibers,
then $Y$ is normal.  (Normality in all finite
characteristics implies normality in characteristic 0).

These strong properties of split varieties
will suffice to prove our Theorem,
provided we construct a compatible splitting.
This is rendered practical by a criterion that was made explicit
 in
\cite{LMP}
in terms of a notion ``residually
normal crossing'', which we now recall.

 A divisor $D$ defined by $f_0=0$ around a point $P$ on
a smooth affine variety $X$ of dimension $n$ has
{\em residually normal crossing} at $P$ if there exists a
 system of parameters $\{x_1,\cdots, x_n\}$ and functions
$f_1,\cdots,f_{n-1}\in k[[x_1,\cdots, x_n]]=\widehat{\OO_P}$ such
 that
$f_i=x_{i+1}f_{i+1}\ ({\rm mod}\ (x_1,\cdots, x_i)
)$ for $i=0,1,\cdots,n-1$, where $f_n=1$ (or a unit).

Now the criterion reads:
\begin{prop}
Let $X$ be a smooth projective variety of dimension $M$ over a
field of characteristic $p > 0$, and let $Z_1, \ldots, Z_M$
be irreducible closed subvarieties of codimension 1.
Suppose that there is a point $P \in X$ such that
$Z_1+ \cdots +  Z_M$
has residually normal crossing at $P$.

{}Further suppose that there exists
a global section s of the anti-canonical bundle $K_X^{-1}$ such
that $\Div s = Z_1 + \cdots + Z_M $.

Then the section $\sigma = s^{p-1}$ gives
a simultaneous compatible
splitting of $Z_1, \ldots,  Z_M$ in $X$.
This  is also a compatible splitting of
any variety obtained from $Z_1,\ldots,Z_M$ by repeatedly
taking intersections and irreducible
components.
\end{prop}

This works because
there is bijection between sections $s$
of $K_X^{-1+p}$
and  $\OO_X$-module morphisms
$\phi: F_* \OO_X \rightarrow \OO_X$,
cf.\ appendix A3 of \cite{vdK}.
In the notation above, one shows that
$\phi$ induces a map $F_* \widehat{\OO_P}/(x_1,\ldots,x_i) \rightarrow
\widehat{\OO_P}/(x_1,\ldots,x_i)$
corresponding with the ``residue'' of $s$ along $x_1=\cdots=x_i=0$
whose divisor
is described by $f_i$.

Several other useful properties of Frobenius splittings can be found
 in
\cite{R2}.

We will apply our theory first in the case $n=3$,
and then indicate the modifications necessary for general $n$.
Instead of directly splitting the pair
$\Ft \subset \Gr(D_3)$, we will find it more convenient
to use an intermediate subspace, which for $n=3$ is just
the triple product of flag varieties $(G/B)^3$.
We embed $\Ftt \subset (G/B)^3$ via the map
$(p_1, p_2, p_3, l_1, l_2, l_3, P)\mapsto ((p_1,l_2),(p_2,l_3),(p_3,l_1))$.
\begin{lem}
The pair $\Ftt \subset (G/B)^3$ is compatibly Frobenius split
\end{lem}
\noindent {\bf Proof.}
We describe the divisor giving our Frobenius splitting on $\Ftt$.
Given $(xB,yB,zB)\in (G/B)^3$, represented by matrices $x$, $y$,
 $z$,
we consider the formula
\begin{eqnarray*}
\llap{$s$}&=&
\left| \begin{array}{ccc}
x_{11} & y_{11} & y_{12} \\
x_{21} & y_{21} & y_{22} \\
x_{31} & y_{31} & y_{32}
\end{array} \right|
\ \cdot \
\left| \begin{array}{ccc}
z_{11} & x_{11} & x_{12} \\
z_{21} & x_{21} & x_{22} \\
z_{31} & x_{31} & x_{32}
\end{array} \right|
\ \cdot \
\left| \begin{array}{ccc}
y_{11} & z_{11} & z_{12} \\
y_{21} & z_{21} & z_{22} \\
y_{31} & z_{31} & z_{32}
\end{array} \right|
\ \cdot \\&&
\left( \ \left| \begin{array}{cc}
x_{21} & x_{22} \\
x_{31} & x_{32}
\end{array} \right|
\cdot
\left| \begin{array}{cc}
y_{11} & y_{12} \\
y_{21} & y_{22}
\end{array} \right|
-
\left| \begin{array}{cc}
x_{11} & x_{12} \\
x_{21} & x_{22}
\end{array} \right|
\cdot
\left| \begin{array}{cc}
y_{21} & y_{22} \\
y_{31} & y_{32}
\end{array} \right|
\ \right)
\ \cdot \\&&
\left| \begin{array}{cc}
x_{11} & y_{11} \\
x_{31} & y_{31}
\end{array} \right|
\ \cdot \
\left| \begin{array}{cc}
z_{11} & z_{12} \\
z_{21} & z_{22}
\end{array} \right|
\ \cdot \
z_{3,1}
\end{eqnarray*}
%
%\begin{eqnarray*}
%\llap{$s$}&=& ( - x_{3,1}y_{1,2}y_{2,1}   +
%     x_{3,1}y_{1,1}y_{2,2} +
%     x_{2,1}y_{1,2}y_{3,1}  -\\&&
%     x_{1,1}y_{2,2}y_{3,1} -
%     x_{2,1}y_{1,1}y_{3,2} +
%     x_{1,1}y_{2,1}y_{3,2} ) \cdot\\&&
%   ( - x_{2,2}x_{3,1}z_{1,1}   +
%     x_{2,1}x_{3,2}z_{1,1} +
%     x_{1,2}x_{3,1}z_{2,1}  \\&&-
%     x_{1,1}x_{3,2}z_{2,1} -
%     x_{1,2}x_{2,1}z_{3,1} +
%     x_{1,1}x_{2,2}z_{3,1} )\cdot \\&&
%   ( - y_{3,1}z_{1,2}z_{2,1}   +
%     y_{3,1}z_{1,1}z_{2,2} +
%     y_{2,1}z_{1,2}z_{3,1} \\&& -
%     y_{1,1}z_{2,2}z_{3,1} -
%     y_{2,1}z_{1,1}z_{3,2} +
%     y_{1,1}z_{2,1}z_{3,2} ) \cdot\\&&
%( x_{3,1}y_{1,1} - x_{1,1}y_{3,1} )\cdot
%   \\&&
%   \left( (-  x_{2,2}x_{3,1}   +
%        x_{2,1}x_{3,2} )   (- y_{1,2}y_{2,1}   +
%        y_{1,1}y_{2,2})\right.-\\&&
%     \left(-  x_{1,2}x_{2,1}   +
%        x_{1,1}x_{2,2} )(  - y_{2,2}y_{3,1}   +
%        y_{2,1}y_{3,2} )\right)\cdot   \\&&
%   ( - z_{1,2}z_{2,1}   +
%     z_{1,1}z_{2,2} )\cdot \\&&z_{3,1}
%\end{eqnarray*}
%
Each factor is a section of some line bundle over $(G/B)^3$.
The first three factors describe the incidence relations $p_1\in
 l_3$,
$p_3\in l_2$, $p_2\in l_1$ respectively.  Their product is a section
of the pull-back to $(G/B)^3$ of the $\OO(1,\dots,1)$ bundle over
$(\PP^2)^3 \times  \Gr(1,\PP^2)^3$. (Just look which Pl\"ucker
coordinates are used.)
Similarly, the last four factors describe a section of that same
 bundle so
that in total $s$ is a section of the anti-canonical bundle of
 $(G/B)^3$,
which is the pullback of $\OO(2,\dots,2)$.

Choosing local coordinates
sensibly, one checks that $\Div s$  has residually normal crossing
 at
the totally degenerate triangle $(E_1,E_1,E_1,E_{12},E_{12},E_{12})$
and we thus have a splitting of $(G/B)^3$ compatible with the components
 of
$\Div s$ and with  the intersection $\Ftt=\II_{3,3}$ of the three
 components
given by the first three factors of $s$.

One also computes that $s$ vanishes
to order five at the locus of totally degenerate triangles.
By \cite{LMP} this implies that the splitting extends to the blowup
of $(G/B)^3$ along that locus, compatibly with the proper transform
 of
$\Ftt$.

\subsection{ Proof of Theorem \ref{vanishing} }

{\bf Proof for $n=3$.}
If $D$ contains
each of the columns of $D_3$ at least once,
then $\LD$ is ample on $(G/B)^3$
and Proposition \ref{splitsur} implies
that $H^i(\FD,\LD)$ vanishes for $i>0$.
{}Furthermore,
$H^0((G/B)^3,\LD)\to H^0(\FD,\LD)$ is surjective,
and a
well-known fact
from representation theory \cite[II 14.20]{J} states
that the restriction map
$H^0((\PP^2)^3 \times  \Gr(1,\PP^2)^3,\LL)\to H^0((G/B)^3,\LL)$
is surjective for any effective $\LL$.
Thus $H^0(\Gr(D), \LD) \rightarrow H^0(\FD,\LD)$ is surjective,
and $H^0(\FD,\LD) \cong S^*_D$ by the previous section,
proving part (a) of the Theorem in this case.

As for part (b), normality follows from the theorem
of Mehta and Srinivas, provided $\Ftt$ possesses
a resolution of singularities with connected fibers.
We will construct two such resolutions in later sections.
The normality, together with the surjectivity of the restriction
 map
above for any multiple of $D$,
is essentially equivalent to projective normality
by \cite{Hart}, Ch II, Ex 5.14(d),
given the projective normality of the Pl\"ucker embedding.

If  $D$ does not contain all the columns of $D_3$, then we need a
strengthening of
Proposition \ref{splitsur}. Indeed Ramanathan has proved \cite[1.12]{R2}
 that
instead of requiring $\LL$ to be ample it suffices
to have a subdivisor $E$ of $\Div s$, not containing any component
 of
 $Y$,
so that $\LL^p\otimes \OO(E)$ is ample. In our case such an $E$ is
provided by the last four factors of $s$, when $\LL=\LD$.

We also need to show that $H^0(\FD,\LD)$ may be identified with
$H^0(\Ftt,\LD)$. For this one considers the projection map
from $(\PP^2)^3 \times  \Gr(1,\PP^2)^3$ to the subproduct
corresponding to the columns that do occur in $D$
(amongst the first six columns
of $D_3$). What one needs to know is that
$\pi_* \OO_{\Ftt} =
% \OO_{\pi\Ftt}$,
\OO_{\FD}$\, ,
where $\pi:\Ftt\to \FD$  is the restriction of the projection map
to $\Ftt=\II_{3,3}$.
Now the fibers of $\pi$
are connected and \cite[6.1.6 and A.1.5]{vdK} apply.

This finishes the proof for $n=3$.
\\[.5em]
\noindent {\bf Case of general $n$.}
\\[.5em]
{}For $n>3$ we need to repeat the argument
``fibered over $\Gr(2,\PP^{n-1})$''.
An element of $\Ft$ may be represented
by a tuple $(M,x,y,z)$ where
 $M$
is an $n \times 3$ matrix of rank three, whose columns span the plane
$P$ in $\PP^{n-1}$, and $x$, $y$, $z$ are three by three matrices
 as
before. The flags $((p_1,l_2)$, $(p_2,l_3)$, $(p_3,l_1))$
are described by the first two columns of
of $Mx$, $My$, $Mz$ respectively.
Let ${P_3}$ denote the stabilizer in $G$ of the plane $E_{123}$.
If $X$ is any ${P_3}$-space, we denote by $G\times^{P_3} X$ the associated
$G$-space fibered over $G/{P_3}=\Gr(2,\PP^{n-1})$,
with fiber $X$ over the point $E_{123}$.

If we replace in the formula for $s$
each $x$ by $Mx$, each $y$ by $My$, each $z$ by $Mz$,
then we get a section of
the relative anti-canonical bundle of the fibration
$$
\begin{array}{ccc}
(GL(3)/B(3))^3 & \rightarrow & G\times^{P_3} (GL(3)/B(3))^3 \\
& & \downarrow \\
& & \Gr(2,\PP^{n-1})
\end{array}
$$
Indeed it has the correct transformation properties under $G$ and
 under
${P_3}$ and it restricts to our known
section of the anti-canonical bundle of the fiber  over $E_{123}$.

Thus to get a section of the anti-canonical bundle of the total space,
one must still
multiply with a section of the pull-back of the anti-canonical
bundle of the base $\Gr(2,\PP^{n-1})$.
There is a choice here. Let us take one for which one can easily
 check
that it has residually normal crossing at $E_{123}$, to wit the
 product
 of
the $n$ Pl\"ucker coordinates (subdeterminants of $M$)
based on taking $3$ consecutive rows of $M$,
with the rows ordered cyclically.

The result of all this is a
section of the anti-canonical bundle of the total space
$G\times^{P_3} (GL(3)/B(3))^3$ that
gives us a splitting with the same virtues as in the case $n=3$.
In particular, it is compatible with $\Ft$.
(It suffices to check this in a
neighborhood of the fiber of $E_{123}$.)

The proof of the Theorem now goes through exactly as before, except
for one problem we still need to address: The analogue of the
``well-known fact from representation theory'' needs to be proved
 now.
We need  to show
that the map
$H^0(\Gr(D_3),\LL)\to H^0(G\times^{P_3} (GL(3)/B(3))^3,\LL)$
is surjective for
effective line bundles on $\Gr(D_3)$.
This is indeed a fact in representation
theory, for which we refer to the Appendix.
Theorem \ref{vanishing} is proved.

\subsection{Fixed points}

In order to gain further information about the
triangle space $\Ft$ and its line bundles, we will
use the method of Lefschetz:  that is, to study the
fixed points of the torus of diagonal
matrices $T \subset GL(n)$ acting on our space.
The work of Atiyah, Bott, and others will then give
us precise formulas for the cohomologies of coherent
sheaves, expressed in terms of the
combinatorial data of the $T$-fixed points and their tangent vectors.
This technique applies only to smooth varieties,
so in subsequent sections
we will study desingularizations of the triangle space.

However, a desingularization map is an isomorphism on the
smooth locus of the variety, so
if a fixed point $\tau$ is a smooth point of $\Ft$,
then the local tangent data will be the same for $\Ft$ and
all desingularizations.  We begin by
examining these smooth points.

We adopt the combinatorial framework for dealing with
general configuration varieties developed in \cite{MaNW},
\cite{MaSchub}.  See also \cite{H}, Lect 16.
A $T$-fixed point of the Grassmannian $\Gr(c,\CC^n)$
is a $c$-dimensional subspace spanned
by coordinate vectors $\{e_{i_1}, \ldots, e_{i_c}\}$,
a subset of
the standard basis $\{e_1,\ldots, e_n\}$ of $\CC^n$.
That is, the fixed points are the spaces $E_I$ for
$I = \{i_1 \leq \cdots \leq i_c\} \subset [1,n]$.
Hence, we may index the fixed points of $\FD$, for
any diagram $D$, by {\em column tabloids} $\tau$,
which are maps $\tau : D \rightarrow [1,n]$,
strictly increasing down each column, such that
for every inclusion of columns $C \subset C'$,
we have $\tau(C) \subset \tau(C')$.
(More precisely, these are the fixed points of $\II_D$,
but in our case this is identical to $\FD$.)
One may check for our $D = D_3$ that there are
$11 n(n-1)(n-2)$ such tabloids.  As we shall see,
all of them are smooth points of $\Ft$ except those
corresponding to maximally degenerate triangles:
namely, the singular fixed points are of the form
$$
% (*) \hspace{1in}
\tau_{ijk} = \left[
\begin{array}{ccccccc}
i &   &    &    & i  & i  & i  \\
  & i &    & i  &    & j  & j  \\
  &   & i  & j  & j  &    & k
\end{array}
\right] \ \ ,
\hspace{1in} \mbox{}
$$
for $i,j,k$ distinct integers in $[1,n]$.

{}For a space $V \in \Gr(c, \CC^n)$, we may model the tangent
space as $T_V \Gr(c,\CC^n) \cong \Hom_{\CC}(V, \CC^n/V )$.
(That is, the tangent bundle of the Grassmannian is isomorphic
to $\Hom$ of the tautological subbundle into the
tautological quotient bundle.)  Furthermore, the incidence
variety
$$
\II = \{(U,V) \in \Gr(c,\CC^n) \times \Gr(c',\CC^n) \mid U \subset
 V \}
$$
has tangent space
$$
T_{(U,V)} \II = \{(\phi,\psi) \in \Hom_{\CC}(U,\CC^n/U)
\times \Hom_{\CC}(V,\CC^n/V) \mid \psi |_U = \phi \mbox{ mod }
 V \} \ .
$$
{}From this, one can deduce as in \cite{MaNW}
\begin{lem}
The eigenvalues of $T$ near a smooth fixed point $\tau$
of $\Ft$
are all of the form $\lam( \diag(x_1,\ldots,x_n) ) = x_i^{-1} x_j$
 for
$i \neq j \in [1,n]$.  The multiplicity $d_{ij}(\tau)$ of
$x_i^{-1} x_j$ is the number of connected components of the
following graph:  vertices = $\{ \mbox{columns } C \mbox{ of }
 D
\mid i\in \tau(C), \ j\not\in \tau(C) \}$,
edges = $\{(C,C') \mid C \subset C' \mbox{ or } C' \subset C \}$.
\end{lem}
{}For instance, consider the singular fixed point $\tau = \tau_{ijk}$
defined above.
Suppose $l,m \neq i,j,k$.  Then the multiplicities
$d_{ab}(\tau)$ are given
by the table:
%$$
%\begin{array}{cccccccc}
%  &   &   &   &   & b &   &    \\
%   & d_{ab} &   & i & j & k & l & m  \\
%   &   &   &-- &-- &-- &-- &--  \\
%   & i & | & 0 & 3 & 1 & 1 & 1  \\
%   & j & | & 0 & 0 & 3 & 1 & 1  \\
% a & k & | & 0 & 0 & 0 & 1 & 1  \\
%   & l & | & 0 & 0 & 0 & 0 & 0  \\
%   & m & | & 0 & 0 & 0 & 0 & 0
%\end{array} \ .
%$$
$$
\begin{array}{r@{\!}l}
\begin{array}{cc}&\\&d_{ab}\end{array}&
\begin{array}{ccccc}&&b & &\\i & j & k & l & m\end{array}\\
\begin{array}{cc}&i\\&j\\a&k\\&l\\&m\end{array}&
\begin{array}{|ccccc}\hline
0 & 3 & 1 & 1 & 1  \\
0 & 0 & 3 & 1 & 1  \\
0 & 0 & 0 & 1 & 1  \\
 0 & 0 & 0 & 0 & 0  \\
 0 & 0 & 0 & 0 & 0
\end{array}
\end{array} \ .
$$
Note that this makes a total of $3n-2$ eigenvectors,
whereas $\Ft$ is $3n-3$ dimensional.  The $3n-2$ vectors
span the Zariski tangent space of $\Ft \subset \Gr(D_3)$
at the singular point.
These eigenvectors also correspond to $T$-stable curves
through $\tau$:  that is,
$\Phi: \CC \rightarrow \Ft$, with $\Phi(0)=\tau$ and
$x\cdot \Phi(s) = \Phi(\lam(x) s)$
for all $s \in \CC$ and $x \in T$, and some eigenvalue character
$\lam$.

{}For all other tabloids, there are $3n-3$ eigenvectors, and the
fixed points are smooth points of $\Ft$.  Now, the singular locus
 of $\Ft$
is $GL(n)$-invariant, and every $GL(n)$-orbit which is not maximally
 degenerate
approaches some $T$-fixed point which is not maximally degenerate,
so we may
conclude:
\begin{lem}
The singular locus of $\Ft$ consists of the maximally
degenerate triangles (for which all three points and all three
lines coincide).
\end{lem}

% \subsection{Picard group}
%
% Using results of M. Brion \cite{Brion}, we can compute the Picard
% group of divisors on $\Ft$, as well as the cones of
% effective and numerically effective divisors in it.
%
%

\section{The Schubert--Semple space}

Next we consider the smooth triangle space $\FSS$ first defined by
Schubert \cite{Schubert}, and given a modern construction
by Semple \cite{Semple}.  See also Collino and Fulton \cite{CF}.
The points of this space may be thought of intuitively
as triangles with an extra piece of data:  a circle passing
through the three vertices, and tangent to the  corresponding
side if two vertices coincide.  The circle is determined by the
triangle in all cases except the maximally degenerate triangles,
for which there is a $\PP^1$ of compatible circles (including
radius zero and infinity).

The rigorous definition is as follows.
{}First, let $n=3$.
The space $Q$ of conic curves
 in $\PP^2$ can be identified as $Q = \PP^*(\Sym^2 \CC^3) \cong
 \PP^5$.
A projective plane in this $\PP^5$ is called a  {\em net}
of conics, and the space of nets is $\Gr(2, Q)$.
Let $\FF^{\circ}_{3,3}$ be the general triangles in $\PP^2$.
{}For any general triangle $\tau$, the conics
passing through its three vertices form a net $N_{\tau}$,
and we have an embedding
$$
\begin{array}{crcl}
\Phi: & \FF^{\circ}_{3,3} & \rightarrow &
\FF^{\circ}_{3,3} \times \Gr(2, Q ) \\
      & \tau       & \mapsto     & (\tau, N_{\tau})
\end{array}.
$$
The Schubert--Semple space is defined as the closure of the
image:
$$
\FF_{3,3}^{\SS} = \mbox{closure } \IM(\Phi) \ \subset \
(\PP^2)^3 \times \Gr(1, \PP^2)^3 \times \Gr(2,Q) \ .
$$

Now, for general $n$, we define $\FSS$ as a family of
such spaces with the plane varying in $\PP^{n-1}$:
$$
\begin{array}{ccc}
\FF_{3,3}^{\SS} & \rightarrow & \FSS \\
& & \downarrow \\
& & \Gr(2,\PP^{n-1})
\end{array}
$$
More formally,
$$
\FSS = GL(n)\times^{P_3}\FF_{3,3}^{\SS}=
{GL(n) \times \FF_{3,3}^{\SS} \over {P_3}} \ .
$$
Here, ${P_3} \subset GL(n)$ is again the parabolic subgroup such
 that
$GL(n)/{P_3} \cong \Gr(2, \PP^{n-1})$, and ${P_3}$ acts on
$GL(n) \times \FF_{3,3}^{\SS}$ by
$p \cdot (g,t) = (gp, b(p^{-1}) t)$, $\ b:{P_3} \rightarrow GL(3)$
being the obvious homomorphism.
Semple shows that this is a smooth projective variety,
and the obvious projection $\FSS \rightarrow \Ft$ is an
isomorphism on the smooth locus of $\Ft$.

We wish to find the $T$-fixed points on $\FSS$, as
well as their tangent eigenvectors.
{}For a point in $\FSS$ whose image is a smooth point of
$\Ft$, this follows immediately from the results of the
previous section.

It remains to consider the fixed points of $\FSS$
lying above the
singular fixed point
$$
\tau_{ijk} = (E_i,E_i,E_i,E_{ij},E_{ij},E_{ij},
E_{ijk}) \in \Ft \ .
$$
The degenerate triangle $\tau_{ijk}$ lies in the plane
$E_{ijk}$, so that
$$
\Sym^2(E_{ijk})^* =
\langle u_i^2, u_j^2, u_k^2, u_i u_j, u_i u_k, u_j u_k \rangle \
 ,
$$
where $u_1, \ldots, u_n$ is the dual of the standard basis of $\CC^n$.
A diagonal matrix $x =\diag(x_1,\ldots,x_n) \in T$ acts on
a monomial by the character
\mbox{$x \cdot u_i u_j = (x_i x_j)^{-1} u_i u_j$}.
According to \cite{CF}, p. 79, the fiber above $\tau_{ijk}$
consists of all nets
$N \in \Gr(2, \PP^*(\Sym^2 E_{ijk}))$
such that
$$
\langle u_k^2, u_j u_k \rangle \subset
N \subset \langle u_k^2, u_j u_k, u_j^2, u_i u_k \rangle \ .
$$
This is a copy of $\PP^1 \subset \FSS$.
More precisely, it is isomorphic as a $T$-space to
$\PP(\CC_{x_j^{-2}} \oplus \CC_{(x_i x_k)^{-1}})$,
where $\CC_{\lam}$ is the one-dimensional representation of $T$
with character $\lam$.
There are exactly two nets in this fiber
which are fixed by $T$, namely
$$\eta_{ijk} = \langle u_k^2, u_j u_k, u_j^2 \rangle = [1:0]
 $$
$$\zeta_{ijk} = \langle u_k^2, u_j u_k, u_i u_k \rangle = [0:1]
 \ . $$
That is, each singular fixed point $\tau_{ijk} \in \Ft$
splits into two isolated fixed points $\eta_{ijk}, \zeta_{ijk}
\in \FSS$.

Now we find the eigenvalues of $T$ acting on the tangent
spaces of these fixed points.
{}First, consider the tangent vector at $\eta_{ijk}$
pointing along the
fiber above $\tau_{ijk}$.
The tangent space of this $\PP^1$ at $\eta_{ijk} = [1:0]$ is
$$
\Hom(\CC_{x_j^{-2}}, \CC_{(x_i x_k)^{-1}})  \cong
\CC_{  x_j^{2} (x_i x_k)^{-1} } \ ,
$$
and the eigenvalue of our tangent vector is
$x_j^{2} (x_i x_k)^{-1}$.
Similarly, $(x_i x_k) x_j^{-2}$ is an eigenvalue of
$T$ acting on the tangent space at $\zeta_{ijk}$.

The other eigenvalues of the smooth tangent spaces
can all be found from examining the
Zariski tangent space of the singular
point $\tau_{ijk}$.  For example, consider
the eigenvector
defined by
\begin{eqnarray*}
\phi \in T_{\tau_{ijk}} \Gr(D_3)& = &\Hom(E_i, \CC^n / E_i)^3 \oplus
\Hom(E_{ij}, \CC^n / E_{ij})^3 \oplus \\
&&\qquad\Hom(E_{ijk}, \CC^n / E_{ijk})^3\ ,
\end{eqnarray*}
\begin{eqnarray*}
\phi& = &(\phi_{ij}, -\phi_{ij}, 0, 0 ,0,0,0) \ ,
\end{eqnarray*}
where $\phi_{ij}(e_l) = \delta_{il} e_j$.
This eigenvector has eigenvalue $x_i^{-1} x_j$
and points along the $T$-stable curve
$\Phi : \CC \rightarrow \Ft$
$$
\Phi(s) = (\Phi_{ij}(s), \Phi_{ij}(-s), E_i, E_{ij}, E_{ij}, E_{ij},
 E_{ijk}) \ ,
$$
where $\Phi_{ij}:\CC \rightarrow \PP^{n-1}$,
$\Phi_{ij}(s) = E_i + s E_j$.
Now, $\Phi(\CC - 0)$ lies in the smooth locus
of $\Ft$, so it lifts uniquely to a $T$-stable curve
$\Phi^{\SS}$ in $\FSS$.
Clearly $\Phi^{\SS}(0)$ is a fixed point above $\tau_{ijk}$,
and we may easily check that it is $\zeta_{ijk}$.
Differentiating $\Phi^{\SS}$ at $s = 0$, we obtain a tangent vector
to $\zeta_{ijk}$ with eigenvalue $x_i^{-1} x_j$.
We may argue similarly for the other Zariski tangent vectors
which do not point along the singular locus of $\Ft$.

{}For a vector which {\em does} point along the singular locus,
for instance
$$
\phi = (\phi_{ij},\phi_{ij},\phi_{ij},0,0,0,0) \ ,
$$
the corresponding curve does not lift uniquely to
a $T$-stable curve in $\FSS$.
In fact, it lifts in exactly two ways, one leading to
each point $\eta_{ijk}$, $\zeta_{ijk}$ above $\tau_{ijk}$.
Hence $\phi$ accounts for a tangent vector with eigenvalue
$x_i^{-1} x_j$ at {\em each} of the lifted fixed points.

Summarizing, we get:
\begin{lem}
\label{SS}
The eigenvalues at the fixed points in $\FSS$
above $\tau_{ijk}$
are as follows:
$$
\begin{array}{cl}
\eta_{ijk} & \ \ \
x_i^{-1} x_k^{-1} x_j^{2}, \
x_i^{-1} x_j, \
x_i^{-1} x_k, \
x_j^{-1} x_k \mbox{ (3 times)}, \
x_i^{-1} x_l, \
x_j^{-1} x_l, \
x_k^{-1} x_l  \\
&  \\
\zeta_{ijk} & \ \ \
x_i x_k x_j^{-2}, \
x_i^{-1} x_j \mbox{ (3 times)}, \
x_i^{-1} x_k, \
x_j^{-1} x_k, \
x_i^{-1} x_l, \
x_j^{-1} x_l, \
x_k^{-1} x_l \ ,
\end{array}
$$
where $l$ runs over $[1,n] \setminus \{i,j,k\}$.
\end{lem}

\section{The Fulton--MacPherson Space}

We describe another desingularization
of $\Ft$, a very special case of
the construction of Fulton and MacPherson in \cite{FM}.
% Locally on $\Ft$ this new resolution $\FFM \rightarrow \Ft$
% is isomorhic to $\FSS \rightarrow \Ft$ as a map, but
%this isomorphism
%involves both source and target of the maps and
%  $\FFM$ comes equipped with a
% different $GL(n)$-action.
It seems this space $\FFM$ is isomorphic to $\FSS$ as an
abstract variety (\cite{FM}, p. 189),
but even if that is so, then it comes equipped with a different
map $\FFM \rightarrow \Ft$ and a different $GL(n)$-action.
In particular, the $T$-fixed points of $\FFM$ are not isolated,
leading to a different type of local data
for our fixed-point formulas below.
Our analysis of the Fulton--MacPherson space
will also show that $\Ft$ has rational singularities.

\subsection{Strata}

Again, we first consider the case $n = 3$, and we will have a
fiber bundle $\FFMt \rightarrow \FFM \rightarrow \Gr(2, \PP^{n-1})$.
It is shown in \cite{FM} how $\FFMt = \PP^2[3]$ can be constructed
as a union of 8 strata,
each consisting of
certain configurations of points and tangent vectors in
$\PP^2$.
{}For each stratum, there is
a natural $GL(3)$ action and an equivariant map to $\Ftt$.
\begin{itemize}
\item
$D_{\emptyset}$, the configuration space
$(\PP^2)^3 \setminus  \bigcup \Delta_{ij}$. \\
An open set in $\FFMt$. \\
Triples $[p_1, p_2, p_3]$ of pairwise-distinct points. \\
A triple maps to the triangle
$$
(p_1,p_2,p_3, \CC p_2 + \CC p_3,
\CC p_1 + \CC p_3, \CC p_1 + \CC p_2) \ ,
$$
where $\CC p + \CC p'$ means the projective line through the points.

\item
$D_{12}$, corresponding to the diagonal
$\Delta_{12} \subset (\PP^2)^3 $. \\
Codimension 1 in $\FFMt$.\\
Configurations
$[p_{12} , p_3$, $\CC^*v_3]$:
distinct points $p_{12}$ and $p_3$,
and a non-zero tangent vector
$v_3\in T_{p_{12}} \PP^{n-1}$, up to
scaling of $v_3$. \\
Intuitively represents infinitesimally distinct points
$
(p_{12}
,\, p_{12} + v_3,\, p_3) \ .
$ \\
Maps to the triangle
$$
(\, p_{12}, p_{12}, p_3,\,
\CC p_{12} + \CC p_3,\, \CC p_{12} + \CC p_3, \,
\CC p_{12} + \CC
v_3\, ) \ .
$$
Similarly $D_{23}$, $D_{13}$.

\item
 $D_{123}$, corresponding to the
total diagonal $\Delta_{123}$ of $(\PP^2)^3$.\\
Codimension 1 in $\FFMt$. \\
Configurations
$[p_{123}, \, \CC^*( v_1, v_2, v_3)]$: a point and three
non-zero tangent
vectors with $v_1 + v_2 + v_3 = 0$, up to
simultaneous scaling.\\
Intuitively represents
$
(p_{123}% + v_1
, p_{123} %+ v_2
+v_3, p_{123} %+ v_3
-v_2) \ .
$ \\
Maps to the triangle
$$
(\, p_{123}, p_{123}, p_{123}, \,
\CC p_{123} + \CC%(v_2 - v_3)
v_1, \,
\CC p_{123} + \CC%(v_1 - v_3)
v_2, \,
\CC p_{123} + \CC%(v_1 - v_2)
v_3 \, ) \ .
$$

\item
$D_{123,12}$. \\
Codimension 2 in $\FFMt$. \\
Configurations
$[p_{123}, \CC^*v_{12},  \CC^*v_3]$:
 a point and two non-zero tangent vectors
up to scaling of each. \\
Intuitively represents
$
(\, p_{123} + v_{12} ,\,
p_{123} + v_{12} + v_3,\,
p_{123} \, ) \ ,
$ with  $v_3$ infinitesimal compared to $v_{12}$, which is
itself already infinitesimal.\\
Maps to  the triangle
$$
(\, p_{123}, p_{123}, p_{123}, \,
\CC p_{123} + \CC v_{12}, \,
\CC p_{123} + \CC v_{12}, \,
\CC p_{123} + \CC v_3 \, ) \ .
$$
Similarly $D_{123,23}$, $D_{123,13}$.

\end{itemize}
The closures of the codimension 1 strata are smooth divisors:
$\overline{D_{ij}} = D_{ij} \cup D_{123,ij}$, \ \
$\overline{D_{123}} = D_{123} \cup D_{123,12} \cup
D_{123,23} \cup D_{123,13}$,
and the intersections are transversal.
The above description is enough to define coordinate charts
for $\FFMt$, gluing together the normal bundles of
the strata appropriately.

Also, the maps described above piece together into
a regular, equivariant desingularization $\pi: \FFMt \rightarrow
 \Ftt$.
The map is an isomorphism outside the singular locus of $\Ftt$,
and the inverse image of a singular triangle
$(p,p,p,l,l,l)$ is the set of configurations:
$[p,  \CC^*(v_1,v_2,v_3)] \in D_{123}$ such that
$v_1$, $v_2$, $v_3$ are non-zero and parallel to $l$,
and $v_1 + v_2 + v_3 = 0$.
(There are also
three extra configurations
in $D_{123,12}$, $D_{123, 23}$, and $D_{123,13}$.)
Thus, the fiber $\pi^{-1}(p,p,p,l,l,l)$ is a projective line,
and it is easily seen that if $(p,p,p,l,l,l)$ is fixed by the torus
 $T$,
then each point of this fiber is also fixed.

\subsection{Fixed lines}

Now let us consider the general $\FFM$, for which everything we have
 said
carries over.  In particular, $\FFM$ has two classes of fixed points:
the isolated ones, which map one-to-one to the
non-singular fixed triangles in
$\Ft$, and the fibers over
the singular fixed triangles $\tau_{ijk}$,
$$
\Pijk = \pi^{-1}(\tau_{ijk}) \ .
$$
The tangent data at the isolated points is identical to that in
$\Ft$.
{}For the $\Pijk$, we will need to determine the types of
their normal bundles
as $T$-equivariant vector bundles over $\PP^1$.

{}First, note that locally near $\Pijk$, we have
$$
\FFM \cong \FFMt(E_{ijk}) \times T_{E_{ijk}} \Gr(2,\PP^{n-1}) .
$$
That is, the normal bundles in the direction of the Grassmannian
are trivial, and we reduce to
$\FFMt(E_{ijk})$, the Fulton--MacPherson space relative to the
plane $\PP^2 = \PP(E_{ijk})$.
Now we will use the
alternative description (\cite{FM}, p. 196)
for $\FFMt$ as a blowup of $(\PP^2)^3$,
first along the triple diagonal $\Delta_{123}$, then along the
proper transforms of the three partial diagonals $\Delta_{12}$,
$\Delta_{23}$,
$\Delta_{13}$.

Let $\CC_{ab}$ denote the one-dimensional $T$-space
on which $\diag(x_1,\ldots,x_n)$ acts by the character $x_a^{-1}
 x_b$.
We may take coordinates for a neighborhood
$U \subset (\PP^2)^3 = \PP(E_{ijk})^3$
near the fixed point $\tau_i = (E_i, E_i, E_i)$
so that, as $T$-spaces, we have
$$
\begin{array}{ccc}
U & \cong &  \left( \CC_{ij}  \times  \CC_{ik} \right)^3  \\
\tau_i & \cong & (0,0,0,0,0,0) \\
\Delta_{12} & \cong & \{ (0,0,a,b,c,d) \} \\
\Delta_{23} & \cong & \{ (a,b,0,0,c,d) \} \\
\Delta_{13} & \cong & \{ (a,b,a,b,c,d) \} \\
\Delta_{123} & \cong & \{ (0,0,0,0,c,d) \} \ .
\end{array}
$$
In this blowup, there will be two $T$-fixed $\PP^1$ above $\tau_i$
 \
(namely, $\tau_{ijk}$ and $\tau_{ikj}$),
and we wish to determine their normal bundles.

Since all the centers of blowing up are products with the last
factor $(\CC_{ij}  \times  \CC_{ik})$, the normal bundles will be
 trivial
in these directions.  Thus, we may reduce to
$$
\begin{array}{ccc}
U' & \cong &  \left( \CC_{ij}  \times  \CC_{ik} \right)^2  \\
\Delta_{12}' & \cong & \{ (0,0,a,b) \} \\
\Delta_{23}' & \cong & \{ (a,b,0,0) \} \\
\Delta_{13}' & \cong & \{ (a,b,a,b) \} \\
\Delta_{123}' & \cong & \{ (0,0,0,0) \} \ .
\end{array}
$$
Performing the first blowup along $\Delta_{123}'$, we obtain:
$$
\begin{array}{ccc}
\Bl_{123} U' & \cong &  \OO(-1) \rightarrow
	       \PP\left( (\CC_{ij}  \times  \CC_{ik})^2 \right)
 \\
\tilde{\Delta}_{12}' & \cong & \OO(-1) \rightarrow
		       \{\, [0:0:a:b] \, \} \\
\tilde{\Delta}_{23}' & \cong & \OO(-1) \rightarrow
		       \{\, [a:b:0:0] \, \} \\
\tilde{\Delta}_{13}' & \cong & \OO(-1) \rightarrow
		       \{ [a:b:a:b]\, \}  \ .
\end{array}
$$
This has two $T$-fixed projective lines: let us focus on
one of them,
$ \PP( \CC_{ij}^2 ) = \{ [a:0:b:0] \}$.  The normal bundle of a
 line $\PP^1
\subset \PP^2$
is $\OO(1)$, so
restricting to a neighborhood $U''$ of our fixed line gives
$$
\begin{array}{ccc}
U'' & \cong &  \OO_{ij}(-1) \oplus 2\OO_{jk}(1) \rightarrow \PP^1
  \\
\Delta_{12}'' & \cong & \{\, (v,0,w) \, \} \rightarrow
		       [0:1]   \\
\Delta_{23}'' & \cong & \{\, (v,w,0) \, \} \rightarrow
		       [1:0] \\
\Delta_{13}'' & \cong & \{\, (v,w,w) \, \} \rightarrow
		       [1:1] \ ,
\end{array}
$$
where $\OO_{ij}(m)$ indicates a line bundle over a $T$-fixed $\PP^1$
with fibers of type $\CC_{ij}$.

Now consider the next blowup, along $\Delta_{12}''$.  This is
locally a product of
$\OO_{ij}(-1) \oplus \OO_{jk}(1)$ and the
locus
$$
(\, 0 \rightarrow [0:1] \, )\  \subset \
(\, \OO_{jk}(1) \rightarrow \PP^1 \, ) \ ,
$$
so the blowup will not affect the first factors, and we may concentrate
on the last.
Thus, consider
the total space of the line bundle $\OO(m)$ over $\PP^1$, and
blow up at a point on the zero-section.
It is easily seen in coordinates that
the normal bundle of the proper transform of the zero-section
is $\OO(m-1)$.

Thus, the second, third, and fourth blowups will transform
$\OO_{ij}(-1) \oplus 2\OO_{jk}(1)$ successively into
$$
\OO_{ij}(-1) \oplus \OO_{jk} \oplus \OO_{jk}(1), \ \
\OO_{ij}(-1) \oplus \OO_{jk} \oplus \OO_{jk}, \ \
\OO_{ij}(-1) \oplus \OO_{jk}(-1) \oplus \OO_{jk}\ .
$$

Recalling the dimensions we dropped at the beginning, we obtain
our final answer.
\begin{lem}\label{bundle}
The normal bundle of the $T$-fixed
component $\Pijk$ in $\FFM$ is
$$
\OO_{ij}(-1) \oplus \OO_{jk}(-1) \oplus \OO_{jk} \oplus
\OO_{ij} \oplus \OO_{jk} \oplus
\sum_{l \in [1,n] \atop l \neq i,j,k}
\left( \OO_{il} \oplus \OO_{jl} \oplus \OO_{kl} \right) \ ,
$$
where $\OO_{ab}(m)$ is the line bundle with Chern class $m$ and
fibers of character $x_a^{-1} x_b$.
\end{lem}

Recall that this describes not just the normal bundle, but an actual
open  neighborhood of the $T$-fixed
component $\Pijk$ in $\FFM$. It follows easily that the canonical
 bundle is
trivial, as a line bundle, on that neighborhood.
% Alternatively, use Grothendieck's theorem on formal functions to
% lift
% a nowhere vanishing section along the fibre to one in a neighborhood.
% And of course one may do it explicitly on the base.
 As $GL(3)$ acts
transitively on the stratum
of degenerate triangles, we get
\begin{lem}\label{free}
The canonical bundle is trivial in a neighborhood of the fiber in
 $\FFM$
of any
singular point.
\end{lem}

\subsection{Rational singularities}

In later sections, it will be convenient to know that
our singular space
$\Ft$  has rational singularities. Let us first recall Kempf's definition
\cite{K}.
A birational proper map $\pi:Y\to X$ is called a {\em rational resolution}\/
if
$Y$ is smooth, and \\
a. $\pi_*\OO_Y=\OO_X$ or, equivalently, $X$ is normal,\\
b. $R^i\pi_*\OO_Y=0$ for $i>0$,\\
c. $R^i\pi_*K_Y=0$ for $i>0$.\\
The last condition is automatic in characteristic $0$.
One says that $X$ has {\em
rational singularities}\/ if there exists a rational resolution $\pi:Y\to
 X$.
The usefulness of this notion lies in the
\begin{lem}
\label{cohom equal}
Let $\pi : Y \to X$ be a  map satisfying conditions
 (a)
and (b),
and $\LL$ a line bundle on $X$.
Then $H^i(X,\LL) = H^i(Y,\pi^*\LL)$ for all $i$.
\end{lem}
This follows from the projection formula and
a degenerate case of the Leray spectral sequence
 \cite[III, Ex. 8.1, 8.3]{H}. We will use the Lemma below
in the case of the triangle space and its desingularizations.

Now, we have seen that the singularity of $\Ft$ is that of the
cone over a quadric in $\PP^3$, and it is well known that this
singularity is rational,
but we shall prove it directly from the definition.

\begin{prop}
The map $\pi:\FFM\to\Ft$ is a rational resolution, and so is
the map $\FSS\to \Ft$.
\end{prop}

\noindent{\bf Proof.}
The target $\Ft$ of $\pi$ is normal by Theorem \ref{vanishing}.
%, for instance because the fibers of
% $\pi$ are connected and $\Ft$ is Frobenius split. (Normality in
% all finite
% characteristics implies normality in characteristic 0).
% As $\pi$ is birational, the first equality follows.
{}For the second condition we use Grothendieck's theorem on formal
 functions.
It tells us that we should try to
show that $H^i(\pi^{-1}(P)_m,\OO/\II^m)$ vanishes, where
 $\pi^{-1}(P)_m$
is the \hbox{$m$-th} order neighborhood of the fiber $\PP^1$
over a point $P$ of the
singular locus and $\II$\/ is the ideal sheaf of this fiber.
By d\'evissage we only need to show that
$H^i(\pi^{-1}(P)_m,\II^{m-1}/\II^m)=H^i(\PP^1,\II^{m-1}/\II^m)$
vanishes for $m>0$. Now $\II^{m-1}/\II^m$ is just a power of the
 conormal
bundle, so by the computation above (lemma \ref{bundle}),
it is a sum of line bundles with
nonnegative Chern class. The result follows.

{}From lemma \ref{free} one sees that the $R^i\pi_*K_\FFM$
are locally the same as the $R^i\pi_*\OO_\FFM$, so they vanish too.
Alternatively, one checks that the Grauert-Riemenschneider vanishing
 theorem
with Frobenius
splitting \cite{MvdK} applies. For this, observe that our splitting
 of
$\Ft$ gives one on the complement of the exceptional locus of $\pi$
 in $\FFM$.
As this exceptional locus has codimension two the splitting extends
 and in fact
our section $s$ of the anti-canonical bundle extends.
The divisor of the extended $s$ contains the proper transform of
 the divisor
of the factor
$
s_4=\left( (-  x_{2,2}x_{3,1}   +
	x_{2,1}x_{3,2} )   (- y_{1,2}y_{2,1}   +
	y_{1,1}y_{2,2})\right.-
     \left(-  x_{1,2}x_{2,1}   +
	x_{1,1}x_{2,2} )(  - y_{2,2}y_{3,1}   +
	y_{2,1}y_{3,2} )\right)
$
of $s$, hence it contains the exceptional locus, as the equation
of the divisor of $s_4$
only puts constraints on two lines in a configuration, no
further restrictions on its points.
\\[.5em]
Before leaving the case of the $\FM$ resolution let us note
that by lemma~\ref{free} there is a line bundle $\omega=\pi_*K_\FFM$
on $\Ft$ whose restriction to the smooth locus is the canonical bundle.
Now let $\pi$ denote the map $\FSS\to \Ft$ instead.
The pull-back of $\omega$ to
$\FSS$ agrees with the canonical bundle outside the exceptional locus,
 which
has codimension two again. It follows that the pull-back is isomorphic
 with
the canonical bundle, so the analogue of lemma \ref{free} holds.
 Thus
the vanishing of $R^i\pi_*K_\FSS$  is equivalent again
to the vanishing of $R^i\pi_*\OO_\FSS$.
To apply the Grauert-Riemenschneider vanishing theorem we now use
 the
factor $s_5=( x_{3,1}y_{1,1} - x_{1,1}y_{3,1} )$, whose divisor
 has a proper
transform containing the exceptional locus.

Another reason that the proposition also holds for the $\SS$
resolution is that it locally looks the same as the $\FM$
resolution.
If locally we
see the singularity as a product of an affine space and a cone over
a product of two projective lines, then it
clearly has an automorphism that interchanges these two lines.
One can pass between the
$\SS$ resolution and the $\FM$ resolution
by means of this local automorphism.
%\\ ******************* \\
%I agree that there is some involution of the singularity which takes
% one
%resolution to the other, but why is it the interchange of the two
% lines?
%Could you (again) send the computations to justify this?
%
%Also, given that the two resolutions ``look the same'' via some
% isomorphism
% of
%the base,
%is it trivial that the rational resolution condition carries over?
%Don't we need that the isomorphism takes the canonical bundle to
% itself?
%Perhaps we {\em should} cut out the discussion of condition (c),
% as you
%suggested,
%and just mention that the singularity is rational in char p as well.
%\\ ******************** \\

\begin{cor}
The desingularizations $\FSS$, $\FFM$, and $\mbox{Bl}_{\mbox{sing
 locus}}
\Ft$ are
all Frobenius split varieties in any characteristic.
\end{cor}

{}Finally, we remark that one can construct the blowup of
$\Ft$ along its singular locus as the
fibered product of $\FFM$ and $\FSS$ over $\Ft$.

\section{Fixed-point formulas}

We apply equivariant fixed-point
theorems to the spaces of the preceding sections,
putting together all the fixed-point data we have accumulated.
This produces explicit formulas
for the $GL(n)$-character and dimension
of the Schur module $S_D$ of any 3-row diagram $D$.
The formulas are more complicated
than those of \cite{MaNW} for northwest diagrams, but
essentially similar.

We discuss the general fixed-point theorems in the first section,
and in the following ones give a summary of the results
in elementary language.  We conclude  by discussing
the possibility of drawing geometric implications from the
combinatorial formulas.

\subsection{General theory}

In what follows, $X$ is
a smooth projective variety
of dimension $M$ over $\CC$,
$L \rightarrow X$ an algebraic line bundle,
and $T = (\CC^*)^n$ a torus acting
on $X$ and $L$.
{\em Throughout this section, we also assume the vanishing
of the higher cohomology groups of $L$}:
$$
H^i(X,L) = 0 \mbox{ for all } i > 0 \ .
$$

The following formula is due to Atiyah and
Bott
\cite{AB}.

\begin{prop}
\label{AB thm}
Suppose the torus $T$
acts $X$ with isolated fixed points.

Then the character of $T$ acting on
the space of global sections of $L$
is given by:
$$
\tr(x \mid H^0(X,L)) =
\sum_{p \ \mbox{\tiny  fixed}} {\tr(x \mid L|_p )
\over  \det( \id -\, x \mid T^*_p X )},
$$
where $p$ runs over the fixed points of $T$, $L|_p$ denotes the
fiber of $L$ above $p$, and $T^*_p X$ is the cotangent space.
\end{prop}

We apply this to
$X = \FSS$ and $L = \pi^* \LD$.
By Lemma \ref{cohom equal} and Theorem~\ref{vanishing},
we have $H^0(\FSS, \pi^* \LD) = H^0(\Ft, \LD) = S^*_D$
and $H^i(\FSS, \pi^* \LD) = H^i(\Ft, \LD) = 0$ for $i>0$,
so the above Proposition gives us a character formula for the Schur
module in terms of the fixed-point data of Lemma \ref{SS}.
We write out the result in the next section.

{}For the $\FM$ space, we need a more general
formula due to Atiyah, Bott, and Singer \cite{AS}.
It requires
the following characteristic
classes for a vector bundle $V$
over any smooth variety $Y$:
$$
{\cal U}^{\lam}(V) =
\prod_{i} {1 - \lam^{-1}\exp(-r_i)
	    \over 1 - \lam^{-1} }
$$
$$
{\cal T}(V) = \prod_{i} {r_i \over 1 - \exp(-r_i)} \ ,
$$
where $\lam$ is a character, and
$r_i$ are the Chern roots of the bundle $V$.
If $V = TY$, we denote ${\cal T}(V) = {\cal T}(Y)$.

\begin{prop}
Suppose ${\cal C}$ is the set of connected
components of the fixed set $X^T$ of $T$.
{}For each component $c \in {\cal C}$,
assume that each restriction
$L|_c$ is a trivial bundle.
Let $N(c) = \bigoplus_{\lam} N_{\lam}(c)$
denote the normal bundle with its
$T$-eigenspace decomposition.

Then the character of $T$ acting on
the space of global sections of $L$
is given by:
$$
\tr(x \mid H^0(X,L)) =
\sum_{c \in {\cal C}}
\left[
{tr(x \mid L|_c) \cdot \prod_{\lam} {\cal U}^{\lam}(N_{\lam}(c))(x)
  \cdot {\cal T}(c)
 \over
 \det( \id -\, x \mid N^*(c) )}
\right](\mbox{\rm Fund }c) \ ,
$$
where the multiplication takes place in the cohomology
ring of the component $c$,
and $\mbox{\rm Fund }c$ denotes the fundamental
homology class.
\end{prop}

Applying this to $X = \FFM$, \, $L = \pi^* \LD$, we
again get a formula for
the character of $S^*_D = H^0(\FFM, \pi^* \LD)$,
this time in terms of the
data in Lemma~\ref{bundle}.

The next result we shall use
is based on the theorem of Hirzebruch-Riemann-Roch
\cite{AS}, combined with Bott's Residue Formula
\cite{Bott}, \cite{AB}, according to the method
of Ellingsrud and Stromme \cite{ES}.
\begin{prop}
Suppose the torus $T = \CC^*$ is one-dimensional,
and acts with isolated fixed points.

Let $\vv = 1$ in the Lie algebra ${\bf t} = \CC$,
and at each $T$-fixed point $p$, let $b(p) = \tr(\vv \mid L|_p
 )$.
Denote the $\vv$-eigenspace decomposition of the
tangent space by $T_p X = \oplus_{i = 1}^M \CC_{r_i(p)} $,
where $r_i(p)$ are the integer eigenvalues.
Also, define the polynomial
$$
\ReRo_M(b;r_1,\dots,r_M) = \mbox{\rm coeff at $U^M \!$ of }
\left( \exp(b U) \prod_{i = 1}^M {r_i U \over 1 -\exp(-r_i U)}
 \right) \ ,
$$
where the right-hand side is considered as a Taylor series in the
formal variable $U$.

Then the dimension of the space of global sections of $L$ is
given by:
$$
\dim  H^0(X,L) =
\sum_{p \ \mbox{\tiny  fixed}}
{\ReRo_M( b(p); r_1(p), \ldots, r_M(p) )
\over  \det( \vv \mid T_p X ) } \ .
$$
\end{prop}

We will consider $X = \FSS$ and take the $T$ in the Proposition
 to be
$\CC^* \subset GL(n)$,\, $q \to \diag(q^{-1},q^{-2},\ldots,q^{-n})
 $.
(This is the principal one-dimensional subtorus corresponding to
 the
half-sum of positive roots.)  Then the eigenvalue characters in Lemma
 \ref{SS}
specialize to the subtorus, and give us the information required
 to
compute the dimension of $S^*_D = H^0(\FSS, \pi^* \LD)$.  (We may
 check
directly that the fixed points of
the subtorus are identical to those of
the large torus of all diagonal matrices.)

Let us also mention that for the smooth spaces $\FSS$ and
$\FFM$, the theorem of Bialynicki-Birula \cite{BB}
gives cell decompositions of these spaces
using the fixed point data.  Thus, one can
compute their singular cohomology groups and Chow groups
as is done in \cite{CF} and \cite{FM}.

\subsection{Character formulas}
\label{character}

{}First, we recall the necessary combinatorial constructions.
We specify a three-row diagram $D$ of squares in the plane
by assigning a multiplicity $m_C \geq 0$ to each column
of the ``universal three-row diagram''
$$
D_3 \  = \ \
\begin{array}{ccccccc}
m_1 & m_2  & m_3 & m_{1,2} & m_{2,3} & m_{1,3} & m_{1,2,3} \\
\Box &       &      & \Box   &        & \Box   & \Box     \\
     & \Box  &      & \Box   & \Box   &        & \Box     \\
     &       & \Box &        & \Box   & \Box   & \Box
\end{array}
$$

We define a {\em standard column tabloid} for $D$\,
with respect to $GL(n)$,
to be a filling (i.e. labeling)
of the squares of $D_3$ by integers in $\{1,\ldots,n\}$,
such that:\\
(i)  the integers in each column are strictly increasing, and\\
(ii) if there is an inclusion  $C \subset  C'$ between two columns,
then all the numbers in the filling of $C$
also appear in the filling of $C'$.
The tabloids describe the fixed points of the
torus $T$ acting on the configuration variety $\Ft$.

Given a tabloid $\tau$ for $D$, define its generating monomial
$$
 x^{\wt(\tau)} = \prod_{(i,j) \in D_3} x_{\tau(i,j)}^{m_j} \ .
$$
That is, the power of $x_i$ is the number of times $i$ appears
in the filling $\tau$, counted with multiplicity.

Also, define integers $d_{ij}(\tau)$ to be the
number of connected components of the following graph:
the vertices are columns
$C$ of $D_3$ such that
$i$ appears in the filling of $C$, but $j$ does
not; the edges are $(C,C')$ such that $C \subset C'$ or $C' \subset
 C$.
\\[1em]
Now, our formula is a sum of terms
corresponding to the column tabloids $\tau$ of $D$.
{}For the smooth tabloids,
the contribution ${\cal C}(\tau)$ is obtained
by the same formula as
in the northwest case discussed in previous works:
$$
{\cal C}(\tau) = {  x^{\wt(\tau)}
\over \prod_{i\neq j} (1-x_i^{-1} x_j)^{d_{ij}(\tau)}  }
$$
However, for the tabloids where
the configuration variety is singular,
we substitute a special contribution
which can be defined in two ways, corresponding
to the two desingularizations $\FSS$ and $\FFM$.
Surprisingly, these expressions reduce algebraically to another,
simpler form which does not appear to be associated with any
desingularization  (c.f. section \ref{virtual}).

\begin{thm}  The character of the Schur module $S_D$ for $GL(n)$
 is
$$
\Char_{S_D} = \sum_{\tau} {\cal C}(\tau) \ ,
$$
where ${\cal C}(\tau)$ are given in the table below.
\end{thm}
To get all tabloids from the types
shown in the table, one should take all permutations of
the first three and the second three columns.
There are $11 n (n-1)(n-2)$ tabloids altogether.

Set
$$
\Delta_{ijk} = \prod_{l \in [1,n] \atop l \neq i,j,k}
\left( 1 - {{x(l)}\over {x(i)}} \right) \,
\left( 1 - {{x(l)}\over {x(j)}} \right) \,
\left( 1 - {{x(l)}\over {x(k)}} \right) \ .
$$

\vfill
\mbox{}
\pagebreak

$$
\begin{array}{ccc}
\mbox{\large \mbox{Tabloid}  $\tau$ }
&
\underline{\mbox{\large \mbox{Character contribution }
${\cal C}(\tau)$ } }
&
\underline{\mbox{\large Desing}}
\\
\overline{
\underline{
\begin{array}{ccccccc}
& & & & & & \\
i &   &   &   & i & i & i \\
  & j &   & j &   & j & j \\
  &   & k & k & k &   & k
\\ & & & & & &
\end{array}
}}
&
{x^{\wt(\tau)} \over
\left( 1 - {{x(i)}\over {x(j)}} \right) \,\left( 1 - {{x(j)}\over
 {x(i)}}
\right) \,
  \left( 1 - {{x(i)}\over {x(k)}} \right) \,\left( 1 - {{x(j)}\over
 {x(k)}}
  \right) \,
  \left( 1 - {{x(k)}\over {x(i)}} \right) \,\left( 1 - {{x(k)}\over
 {x(j)}}
  \right)
\Delta_{ijk}
}
& \mbox{smooth}
\\
\underline{
\begin{array}{ccccccc}
& & & & & & \\
i &   &   &   & i & i & i \\
  & i &   & i &   & j & j \\
  &   & j & j & j &   & k
\\ & & & & & &
\end{array}
}
&
{x^{\wt(\tau)} \over
\left( 1 - {{x(i)}\over {x(j)}} \right) \,\left( 1 - {{x(j)}\over
 {x(i)}}
\right) \,
  \left( 1 - {{x(j)}\over {x(k)}} \right) \,{{\left( 1 - {{x(k)}\over
 {x(i)}}
  \right) }^2}\,
  \left( 1 - {{x(k)}\over {x(j)}} \right)
\Delta_{ijk}
}
& \mbox{smooth}
\\
\underline{
\begin{array}{ccccccc}
& & & & & & \\
i &   &   &   & i & i & i \\
  & i &   & i &   & k & j \\
  &   & j & j & j &   & k
\\ & & & & & &
\end{array}
}
&
{x^{\wt(\tau)} \over
\left( 1 - {{x(i)}\over {x(j)}} \right) \,{{\left( 1 - {{x(j)}\over
 {x(i)}}
\right) }^2}\,
  \left( 1 - {{x(k)}\over {x(i)}} \right) \,{{\left( 1 - {{x(k)}\over
 {x(j)}}
  \right) }^2}
\Delta_{ijk}
}
& \mbox{smooth}
\\
\underline{
\begin{array}{ccccccc}
& & & & & & \\
i &   &   &   & i & i & i \\
  & i &   & i &   & k & j \\
  &   & i & j & j &   & k
\\ & & & & & &
\end{array}
}
&
{x^{\wt(\tau)} \over
{{\left( 1 - {{x(j)}\over {x(i)}} \right) }^2}\,\left( 1 - {{x(j)}\over
 {x(k)}}
\right) \,
  \left( 1 - {{x(k)}\over {x(i)}} \right) \,{{\left( 1 - {{x(k)}\over
 {x(j)}}
  \right) }^2}
\Delta_{ijk}
}
& \mbox{smooth}
\\
\underline{
\begin{array}{ccccccc}
& & & & & & \\
i &   &   &   & i & i & i \\
  & i &   & i &   & j & j \\
  &   & i & j & j &   & k
\\ & & & & & &
\end{array}
}
&
{x^{\wt(\tau)}
\left( 1 - \left( { 1 - {x(i)\over x(j)}  }\right)^{-1}
  - \left( { 1 - {x(j)\over x(k)}  }\right)^{-1}
\right)
\over
{{\left( 1 - {{x(j)}\over {x(i)}} \right) }^2}\,\left( 1 - {{x(k)}\over
 {x(i)}}
\right) \,
  {{\left( 1 - {{x(k)}\over {x(j)}} \right) }^2}
\Delta_{ijk}
}
& \FM
\\
 &
=
{x^{\wt(\tau)} \over
 {\left( 1 - {x(j)\over x(i)} \right)}^3\,
  \left( 1 - {x(k)\over x(i)} \right)   \,
  \left( 1 - {x(k)\over x(j)} \right)   \,
  \left( 1 - {x(i)\,x(k) \over x(j)^2} \right)
\Delta_{ijk}
}
\ \ +
& \SS
\\
& \ \ \ \ \
{x^{\wt(\tau)} \over
  \left( 1 - {x(j)\over x(i)} \right)   \,
  \left( 1 - {x(k)\over x(i)} \right)   \,
 {\left( 1 - {x(k)\over x(j)} \right)}^3\,
  \left( 1 - {x(j)^2 \over x(i)\,x(k)} \right)
\Delta_{ijk}
}
\\[3em]
 &
=
{x^{\wt(\tau)} \over
 {\left( 1 - {x(j)\over x(i)} \right)}^3\,
 {\left( 1 - {x(k)\over x(j)} \right)}^3\,
\Delta_{ijk}
}
& \mbox{\bf ?}
\end{array}
$$

\pagebreak
\normalsize

\subsection{Dimension formula}
\label{dimension}

We compute the dimension of the $GL(n)$-module $S_D$.
Again, each tabloid gives a contribution, which is
the value of a certain multivariable polynomial
$\ReRo$ at a sequence
of integers specific to the tabloid.

Let $M = 3n-3$, the dimension of the triangle space $\Ft$.
Define an $(M+1)$-variable polynomial, homogeneous of degree $M$,
 by
$$
\ReRo_M(b;r_1,\dots,r_M) = \mbox{coeff at $U^M \!$ of }
\left( \exp(bU) \prod_{i = 1}^M {r_i U \over 1 - \exp(-r_i U)}
 \right) \ ,
$$
where the right side is understood as a Taylor series in $U$.
{}For example, for $n=3$, $M=6$, $\ReRo_6(b;r_1,\ldots,r_6)$
is an irreducible 7-variable polynomial,
homogeneous of degree 6, with 567 terms.
 However, since we will only
evaluate this polynomial at $(M+1)$-tuples
of small integers, this is
within the range of computer
calculations provided the column multiplicities
$m_j$ are not very large.

{}For the smooth tabloids, we again have a formula
for the contributions
in terms of the integers $d_{ij}(\tau)$.
Namely, define a multiset of integers
$r(\tau) = \{r_1, r_2, \ldots, r_{M}\}$ by inserting
the entry $i-j$ with multiplicity $d_{ij}(\tau)$
for each ordered pair $i, j \in [1,n]$.
That is, the total multiplicity of the integer $k$ in $r(\tau)$ is
$$
 \sum_{i,j \in [1,n] \atop i-j = k} d_{ij}(\tau) \ .
$$
Let $b(\tau)$ be the sum of the entries of the tabloid $\tau$, counting
column multiplicities:
$$
b(\tau) = \sum_{(i,j) \in D} m_j\, \tau(i,j) \ .
$$
Then the contribution of $\tau$ to the dimension formula is
\begin{eqnarray*}
{\cal R}(\tau) & = & {\ReRo_M(b(\tau); r(\tau)) \over \prod_{k
 = 1}^M
r(\tau)_i} \\[.2em]
& = & {\ReRo_M(b(\tau); r(\tau)) \over \prod_{i,j \in [1,n]} (i-j)}
\end{eqnarray*}

{}For the singular tabloids, we have only
an expression corresponding to the $\SS$ desingularization,
as well as a simplified expression with no geometric explanation,
as before.
% (The authors have verified
% this simplified expression by computation for $n = 3,4$,
% but have not proved it for general $n$.)

\begin{thm}  The dimension of the Schur module $S_D$ for $GL(n)$
 is
$$
\dim S_D = \sum_{\tau} {\cal R}(\tau) \ ,
$$
where ${\cal R}(\tau)$ are given in the table below.
\end{thm}

The entries in the table are terms of the form
${\cal R} = \ReRo_M(b;r) / \prod r$ for integers $b$ and multisets
 $r$.
Let $r'(i,j,k)$ be the standard multiset
with entries $i-l$, $j-l$, $k-l$ for each integer $l \in [1,n]$,
 $l \neq i,j,k$.
(For example, if $n=4$, we have $r'(1,2,3) = \{ 1-4,\, 2-4,\,
 3-4 \} =
\{-3,-2,-1\}$.)

\vfill
\mbox{}
\pagebreak

$$
\large
\begin{array}{c@{\!}c@{\!}c}
\mbox{Tabloid } \tau &
\underline{\mbox{Dimension contribution }{\cal R}(\tau) } &
\underline{\mbox{Desing}}\\
\overline{
\underline{
\begin{array}{ccccccc}
& & & & & & \\
i &   &   &   & i & i & i \\
  & j &   & j &   & j & j \\
  &   & k & k & k &   & k
\\ & & & & & &
\end{array}
}}
&
{\ReRo_M(b(\tau);\ j-i,\ i-j,\ k-i,\ k-j,\ i-k,\ j-k,\ r'(i,j,k)
 )
\over (j-i)(i-j)(k-i)(k-j)(i-k)(j-k) \ \prod r'(i,j,k) }
& \mbox{smooth}
\\
\underline{
\begin{array}{ccccccc}
& & & & & & \\
i &   &   &   & i & i & i \\
  & i &   & i &   & j & j \\
  &   & j & j & j &   & k
\\ & & & & & &
\end{array}
}
&
{\ReRo_M(b(\tau);\ j-i,\ i-j,\ k-j,\ i-k,\ i-k,\ j-k,\ r'(i,j,k)
 )
\over (j-i)(i-j)(k-j)(i-k)^2 (j-k) \ \prod r'(i,j,k) }
& \mbox{smooth}
\\
\underline{
\begin{array}{ccccccc}
& & & & & & \\
i &   &   &   & i & i & i \\
  & i &   & i &   & k & j \\
  &   & j & j & j &   & k
\\ & & & & & &
\end{array}
}
&
{\ReRo_M(b(\tau);\ j-i,\ i-j,\ i-j,\ i-k,\ j-k,\ j-k,\ r'(i,j,k)
 )
\over (j-i)(i-j)^2 (i-k) (j-k)^2 \ \prod r'(i,j,k) }
& \mbox{smooth}
\\
\underline{
\begin{array}{ccccccc}
& & & & & & \\
i &   &   &   & i & i & i \\
  & i &   & i &   & k & j \\
  &   & i & j & j &   & k
\\ & & & & & &
\end{array}
}
&
{\ReRo_M(b(\tau);\ i-j,\ i-j,\ k-j,\ i-k,\ j-k,\ j-k,\ r'(i,j,k)
 )
\over (i-j)^2 (k-j)(i-k)(j-k)^2 \ \prod r'(i,j,k) }
& \mbox{smooth}
\\
\underline{
\begin{array}{ccccccc}
& & & & & & \\
i &   &   &   & i & i & i \\
  & i &   & i &   & j & j \\
  &   & i & j & j &   & k
\\ & & & & & &
\end{array}
}
&
\begin{array}{r}
{\ReRo_M(b(\tau);\ i-j,\ i-j,\ i-j,\ i-k,\ j-k,\ 2j-i-k,\ r'(i,j,k)
 )
\over (i-j)^3 (i-k)(j-k)(2j-i-k) \ \prod r'(i,j,k) } \  +  \ \\[.4em]
{\ReRo_M(b(\tau);\ i-j,\ i-k,\ j-k,\ j-k,\ j-k,\ i+k-2j,\ r'(i,j,k)
 )
\over (i-j)(i-k)(j-k)^3 (i+k-2j) \ \prod r'(i,j,k) }
\end{array}
& \SS
\\[1em]
 &
=  {\ReRo_M(b(\tau);\ i-j,\ i-j,\ i-j,\ j-k,\ j-k,\ j-k,\ r'(i,j,k)
 )
\over (i-j)^3 (j-k)^3 \ \prod r'(i,j,k) }
& \mbox{\bf ?}
\end{array}
$$

% \pagebreak
\normalsize

\vfill \mbox{}
\pagebreak

\subsection{Virtual desingularization}
\label{virtual}

In the previous sections, we have drawn consequences about Schur
 modules
from the geometry of the triangle space.  However, one can imagine
 reversing
this process.

{}For instance, consider the trivial line bundle over the triangle
space, $m_1 = m_2 = \cdots = m_7 = 0$.  Its space of sections
 is the
trivial Schur module,
with character 1.  Hence, our character formulas for the $\FM$ and
 $\SS$
desingularizations each give
ludicrously complicated expressions for 1.
There is a non-trivial contribution
${\cal C}(\tau)$ for each $T$-fixed point $\tau$ of the smooth space:
the numerators are reduced to 1, but the denominators still possess
 a factor for each eigenvalue of the tangent space at the fixed point.
Because these eigenvalues determine a cell decomposition of the smooth
space, one can read off from these ``ludicrous formulas''
a great deal of geometric information about the desingularizations
 of $\Ft$.

In fact, suppose we did not know of the existence of the $\SS$
desingularization.
We could guess that there exists such a smooth space, with its cell
 decomposition,
by finding a ludicrous formula for 1 with terms of Atiyah-Bott type.

This is not difficult:
let us start by assuming there exists a $GL(n)$-equivariant
desingularization of $\Ftt$
with fibers of dimension 1 (the smallest possible), and with each
 singular
fixed point lifting to only two fixed points in the smooth space.
The singular fibers must then be isomorphic to $\PP^1$, since this
 is the
only curve
possessing an appropriate torus action.
Now, the contributions
${\cal C}(\tau)$ from the
smooth tabloids of $\Ft$ are constrained to be the entries in the
 table.
As for the singular tabloids, they must each correspond to two terms
 in the
formula,
one for each lifted fixed point.  We know most of their eigenvalues
 from the
$T$-stable curves in $\Ft$: each will lift to either
two or one $T$-stable curves in the smooth space, depending on whether
or not it lies in the singular locus.  Since there are six eigenvectors
at each fixed point of the smooth space, this leaves only one eigenvector
to determine at each of the two lifted fixed points.
Since the fiber is $\PP^1$, these must be reciprocals of each other.
Now, the six singular tabloids lie in a single GL(n) orbit, so the
corresponding unknown eigenvalues are all images of each other.
Therefore, there is only one variable left unknown, which we can
 solve
for in the ludicrous equation:  over $\tau_{123}$ the value
must be $x_1 x_3 x_2^{-2}$, as in our formula.

Now, for the space of tetrahedra and higher-dimensional simplices
(corresponding to general diagrams with four or more rows), there
is no known explicit desingularization.  One may hope to find evidence
of one by arguments like the above, combined with induction on the
 dimension.
That is, an appropriate ludicrous formula for 1 may be considered
 as a
combinatorial
or virtual desingularization.

{}Finally, let us point out one mystery in our results:  the algebraic
simplification
of the character and dimension formulas, combining the complicated
 contributions
given by our desingularizations into a single term of Atiyah-Bott
 type.
In the above philosophy, this would mean that $\Ft$ is already ``virtually
smooth'',
with eigenvalues above $\tau_{ijk}$ equal to $x_i^{-1} x_j$ (three
 times), and
$x_j^{-1} x_k$ (three times).  Note that there can exist no actual
$G$-equivariant desingularization
of $\Ft$ for which the singular fixed points  each lift uniquely,
 since
this would
use up only six of the seven eigenvalues given by the $T$-stable
curves at $\tau_{ijk}$ (cf. Lemma \ref{SS}).
(Each of these curves must lift to at least one $T$-stable curve
 in the
smooth space
% (since the curves are translates of the fixed point by
% a T-stable curve in G itself)
and lead to some lifted fixed point.)

\section{Appendix: Restriction and induction}\label{resind}
We now prove a result from the representation theory of reductive
 algebraic
groups, needed in the proof of theorem \ref{vanishing}.
Let $G$ be a connected reductive algebraic
group, $B$ a Borel subgroup, $P$ a parabolic subgroup
containing $B$.
Let us call a weight $\lambda$ {\em effective} if $\ind_B^G\lambda\neq0$.
(In \cite{J} an effective weight is called dominant, in \cite{vdK}
it is called anti-dominant.) For the notions of `induction', `good
 filtration',
`excellent filtration', and the basic theorems concerning them we
 refer to
\cite{J}, \cite{Mat}, \cite{vdK}.

\begin{lem}
Let $\lambda$ be effective and let $M$ be a $P$-module that has excellent
filtration (as a $B$-module).
Then $\ind_B^P(\lambda)\otimes M$ has an excellent filtration.
\end{lem}

\noindent {\bf Proof.} As
$\ind_B^P(\lambda)\otimes M=\ind_B^P(\lambda\otimes M)$
this follows from the main theorems
on excellent filtrations as collected
in \cite{vdK}.

\begin{prop}
{}For effective $\lambda_1$, \ldots, $\lambda_n$, the restriction
 map
$$
\res :\ind_B^G(\lambda_1) \otimes \cdots \otimes \ind_B^G(\lambda_n)
 \to
\ind_B^G(\ind_B^P(\lambda_1) \otimes \cdots \otimes \ind_B^P(\lambda_n))
$$
is surjective.
\end{prop}

\noindent {\bf Proof.}
It suffices to show that for each $i$ the kernel $\ker\phi_i$ of
the surjective map $\phi_i$:
$$\begin{array}{c}
\ind_B^G(\lambda_1)\otimes\cdots\otimes\ind_B^G(\lambda_{i-1})
\otimes\ind_B^G(\lambda_{i})\otimes\ind_B^P(\lambda_{i+1})
\otimes\cdots\otimes
\ind_B^P(\lambda_n)\\
\downarrow\\
\ind_B^G(\lambda_1)\otimes\cdots\otimes\ind_B^G(\lambda_{i-1})
\otimes\ind_B^P(\lambda_{i})\otimes\ind_B^P(\lambda_{i+1})
\otimes\cdots\otimes
\ind_B^P(\lambda_n)
\end{array}$$
is $\ind_B^G$-acyclic. Indeed $\res$ may be viewed as
$\ind_B^G(\phi_1)\circ\cdots\circ\ind_B^G(\phi_n)$. Now a module
 $M$ is
$\ind_B^G$-acyclic
if and only if $k[G]\otimes M$ is $B$-acyclic, so that the result
 follows
from the lemma and the
main theorems on excellent filtrations {\it etc.}

\end{document}